\documentclass{gtart}

\usepackage{amsmath,amssymb}
\usepackage{graphicx}

\newtheorem{theorem}{Theorem}[section]
\newtheorem{lemma}[theorem]{Lemma}

\theoremstyle{definition}

\theoremstyle{remark}

\numberwithin{equation}{section}

\begin{document}

\title{On parameterizations of Teichm\"uller spaces of surfaces with boundary}

\author{Ren Guo}

\address{Department of Mathematics, Rutgers University, Piscataway, NJ, 08854, USA}

\email{renguo@math.rutgers.edu}

\begin{abstract}

In \cite{rigidity}, Luo introduced a $\psi_{\lambda}$ edge
invariant which turns out to be a coordinate of the Teichm\"uller
space of a surface with boundary. And he proved that for $\lambda
\geq 0$, the image of the Teichm\"uller space under
$\psi_{\lambda}$ edge invariant coordinate is an open cell. In
this paper we verify his conjecture that for $\lambda<0$, the
image of the Teichm\"uller space is a bounded convex polytope.

\end{abstract}


\maketitle

\section{Introduction}

Suppose $S$ is a compact connected surface of non-empty boundary
and has negative Euler characteristic. It is well known that there
are hyperbolic metrics with totally geodesic boundary on the
surface $S$. Two such hyperbolic metrics are called \it isotopic
\rm if there is an isometry isotopic to the identity between them.
The space of all isotopy classes of hyperbolic metrics on $S$,
denoted by $T(S)$, is called the Teichm\"uller space of the
surface $S$.

There are several known parameterizations of the Teichm\"uller
spaces. In particular, using the 3-holed decomposition of a
surface, Fenchel-Nielsen introduced a coordinate for $T(S)$, for
more detail see the book Imayoshi $\&$ Taniguchi \cite{it}.
Bonahon \cite{bo} produced a parametrization of the Teichm\"uller
spaces using the sheared coordinate. Penner \cite{p1,p2}
introduced the ``lambda length" coordinate and simplicial
coordinate of the decorated Teichm\"uller space. Recently Luo
\cite{boundary,rigidity} introduced a family of coordinates of
$T(S)$. To be more precise, for each real number $\lambda$, he
introduced a $\psi_{\lambda}$ edge invariant associated to a
hyperbolic metric which turns out to be a coordinate of the
Teichm\"uller space $T(S)$. When $\lambda \geq 0$, he proved that
the image of the Teichm\"uller space under the coordinate is an
open convex polytope independent of $\lambda.$ Luo \cite{luo}
conjectured that for $\lambda<0$, the image of the Teichm\"uller
space under $\psi_{\lambda}$ edge invariant coordinate is a
bounded convex polytope. The purpose of this paper is to verify
this conjecture.

Let us begin by recaling the $\psi_{\lambda}$ edge invariant
coordinate introduced by Luo \cite{rigidity}. The coordinate
depends on a fixed ideal triangulation of $S.$ Recall that a
colored hexagon is a hexagon with three non-pairwise adjacent
edges labelled by red and the opposite edges labelled by black.
Take a finite disjoint union of colored hexagons and identify all
red edges in pairs by homeomorphisms. The quotient is a compact
surface with non-empty boundary together with an ideal
triangulation. The 2-cells in the ideal triangulation are
quotients of the hexagons. The quotients of red edges
(respectively black edges) are called the edges (respectively \it
A-arcs\rm) of the ideal triangulation. It it well known that every
compact surface $S$ of non-empty boundary and negative Euler
characteristic admits an ideal triangulation.

In a hyperbolic metric, any hexagon in an ideal triangulation is
isotopic (leaving the boundary of a surface fixed) to a hyperbolic
right-angled hexagon. It is well known that a hyperbolic
right-angled hexagon is determined up to isometry preserving
coloring by the lengths of three red edges. Furthermore, for any
$l_1, l_2, l_3 \in \mathbf{R}_{>0}$, there exists a unique colored
hyperbolic right-angled hexagon whose three red edges have lengths
$l_1, l_2, l_3$, for a proof see Buser \cite{bu}.

Given an ideally triangulated surface $S$ with $E$ the set of all
edges, each hyperbolic metric $d$ on $S$ has a length coordinate
$l_d: E \to \mathbf{R}_{>0}$ which assigns each edge $e$ the
length of the shortest geodesic arc homotopic to $e$ relative to
the boundary of $S$.  On the other hand, given a function $l: E
\to \mathbf{R}_{>0}$, we can produce a hyperbolic metric with
totally geodesic boundary on $S$. This metric is constructed by
making each 2-cell with red edges $e_i, e_j, e_k$ a colored
hyperbolic right-angled hexagon the lengths of whose red edges are
$l(e_i), l(e_j), l(e_k)$. Thus, the Teichm\"uller space $T(S)$ can
be identified with the space $\mathbf{R}^E_{>0}$ by length
coordinates.

\begin{figure}[htbp]
\begin{center}
\includegraphics[scale=.4]{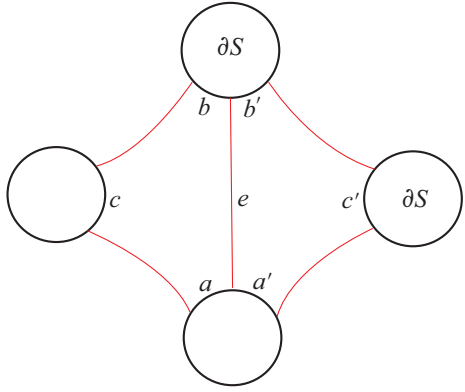}
\end{center}
\caption{\label{invariant}The lengths of A-arcs in the definition of
$\psi_{\lambda}$ edge invariant are labeled.}
\end{figure}

In \cite{rigidity}, Luo introduced the \it $\psi_{\lambda}$ edge
invariant \rm of a hyperbolic metric as $\psi_{\lambda}: E
\rightarrow \mathbf{R}$ defined by
$$\psi_\lambda(e)=\int_0^{\frac{a+b-c}{2}} \cosh^{\lambda}(t) dt +
\int_0^{\frac{a'+b'-c'}{2}} \cosh^{\lambda}(t)dt$$ where $e$ is an
edge of an ideal triangulation shared by two hyperbolic
right-angled hexagons and $a,b,c,a',b',c'$ are lengths of the
A-arcs labelled as in Figure \ref{invariant}. Now consider the map
$\Psi_{\lambda}: T(S) \to \mathbf{R}^E$ sending a hyperbolic
metric $l$ to its $\psi_{\lambda}$ edge invariant.

The following two theorems are proved in Luo \cite{rigidity}. The
special case of $\lambda=0$ was proved in Luo \cite{boundary}. We
use $(S,T)$ to denote a surface $S$ with an ideal triangulation
$T.$
\begin{theorem}(Luo)\label{embedding}
Suppose $(S, T)$ is an ideally triangulated surface. For any
$\lambda \in  \mathbf{R}$, the map $\Psi_{\lambda}: T(S) \to
\mathbf{R}^E$ is a smooth embedding. In particular, each
hyperbolic metric with geodesic boundary on $(S, T)$ is determined
up to triangulation preserving isometry by its $\psi_{\lambda}$
edge invariant.
\end{theorem}

An \it edge path \rm $(H_0, e_1, H_1,..., e_n, H_n)$ is a
collection of hexagons and edges in an ideal triangulation so that
two distinct hexagons $H_{i-1}$ and $H_i$ sharing the edge $e_i$
for $i=1,...,n.$ An edge path $(H_0, e_1, H_1,..., e_n, H_n)$ is
an \it edge cycle \rm if $H_0=H_n$. For example see Figure
\ref{path}. A \it fundamental edge path (or fundament edge cycle)
\rm is an edge path (or edge cycle) so that each edge in the ideal
triangulation appears at most twice in the path (or cycle).

\begin{theorem}(Luo)\label{positive} Let $\lambda \geq 0$. For an
ideal triangulated surface $(S, T)$, $ \Psi_{\lambda}(T(S)) =\{ z
\in \mathbf{R}^E | $ for each fundamental edge cycle $(H_0, e_1,
H_1,..., e_n, H_n=H_0)$, $\sum_{i=1}^n z(e_i)
>0$\}. Thus $\Psi_{\lambda}(T(S))$ is an open convex polytope independent of the
parameter $\lambda \geq 0$.
\end{theorem}

In this paper we generalize Theorem \ref{positive} to any real
number $\lambda.$ The main result is the following.

\begin{theorem}\label{negative}
For an ideal triangulated surface $(S, T)$, $
\Psi_{\lambda}(T(S))$ is the set of points $z \in \mathbf{R}^E$
satisfying

1. $z(e) < 2 \int_{0}^{\infty} \cosh^{\lambda}(t) dt$ for each
edge e;

2. $\sum_{i=1}^n z(e_i) >-2 \int_{0}^{\infty} \cosh^{\lambda}(t)
dt$ for each fundamental edge path $(H_0, e_1, H_1,...,$ $e_n,
H_n)$;

3. $\sum_{i=1}^n z(e_i) >0$ for each fundamental edge cycle $(H_0,
e_1, H_1,..., e_n, H_n=H_0)$.

Thus $\Psi_{\lambda}(T(S))$ is an open convex polytope. And
$\Psi_{\lambda_1}(T(S))\subset
\Psi_{\lambda_2}(T(S))\subset\Psi_0(T(S))=\Psi_{\lambda_3}(T(S))$
for $\lambda_1<\lambda_2<0<\lambda_3.$ The intersection
$\cap_{\lambda=0}^{-\infty}\Psi_{\lambda}(T(S))$ is empty.
\end{theorem}

It is easy to see when $\lambda\geq 0$ the conditions in Theorem
\ref{negative} are reduced to the third one which is exactly the
condition in Theorem \ref{positive}. The proof of Theorem
\ref{negative} follows the same strategy used in Luo's proof of
Theorem \ref{positive} \cite{rigidity}.

In section 2 we investigate degenerations of a hyperbolic
right-angled hexagon. In section 3 we prove the main result
Theorem \ref{negative}.

\section{Degenerations of a hyperbolic hexagon}

\begin{figure}[htbp]
\begin{center}
\includegraphics[scale=.35]{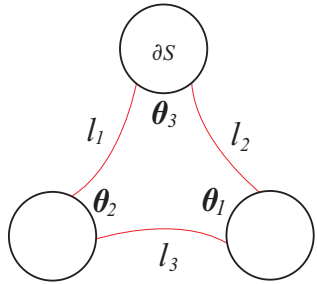}
\end{center}
\caption{\label{ltheta}An hyperbolic
right-angled hexagon with lengths of edges and A-arcs labeled.}
\end{figure}

In this section we always assume a hyperbolic right-angled hexagon
has three non-pairwise adjacent edges of lengths $l_1, l_2, l_3$
and opposite A-arcs of lengths $\theta_1, \theta_2, \theta_3$
labelled in Figure \ref{ltheta}. And recall that the r-coordinate
is defined as $r_i=\frac{\theta_j+\theta_k-\theta_i}{2}.$

We improve a lemma proved in Luo \cite{rigidity}.
\begin{lemma}\label{r} Consider $r_i$ as a function of $(l_1,l_2,l_3).$ We have
$\lim_{l_i \to 0} r_i =\infty$ so that the convergence is uniform
in $(l_1,l_2,l_3)$.
\end{lemma}

\begin{proof}  By the cosine law of a hyperbolic right-angled hexagon, we
see that for $i\neq j\neq k\neq i,$
\begin{eqnarray*}
\cosh\theta_j &=& \frac{ \cosh l_j +
\cosh l_i\cosh l_k}{\sinh l_i \sinh l_k}\\
& > & \frac{ \cosh l_i
\cosh l_k}{\sinh l_i \sinh l_k }\\
& \geq & \frac{ \cosh l_i }{\sinh l_i }.
\end{eqnarray*}
Hence we have $\lim_{l_i \to 0} \theta_j =\infty.$ Thus $\lim_{l_i
\to 0} \frac{\cosh\theta_j}{\sinh\theta_j}=1.$ By symmetry we have
$\lim_{l_i \to 0} \frac{\cosh\theta_k}{\sinh\theta_k}=1.$

On the other hand, by cosine law we see that for $i\neq j\neq
k\neq i,$
$$\cosh l_i-\frac{\cosh \theta_j\cosh \theta_k}{\sinh \theta_j \sinh \theta_k}=
\frac{\cosh \theta_i}{\sinh \theta_j \sinh \theta_k}>
\frac{2e^{\theta_i}}{e^{\theta_j+\theta_k}}=\frac2{e^{2r_i}}.$$
Since the left hand side converges to 0 as $l_i\to 0,$ we have
$\lim_{l_i \to 0} r_i =\infty$.

To show the convergence is uniform, we consider the following
formula called \it tangent law \rm derived in Luo \cite{rigidity}.
For $i\neq j\neq k\neq i,$
$$\tanh^2\frac{l_i}2 =
\frac{\cosh r_j\cosh r_k}{\cosh r_i\cosh (r_i+r_j+r_k)}.$$ By the
formula,

\begin{eqnarray*}
\tanh^2\frac{l_i}2 &=& \frac{1}{\cosh r_i}\cdot\frac1{(1+\tanh r_j
\tanh r_k)\cosh r_i
+(\tanh r_j  +\tanh r_k)\sinh r_i }\\
&\geq& \frac{1}{\cosh r_i}\cdot\frac1{(1+1)\cosh r_i + (1+1)|\sinh r_i|}\\
&\geq& \frac{1}{4 \cosh^2 r_i}.
\end{eqnarray*}

It follows that $$\cosh^2 r_i \geq \frac{1}{4
\tanh^2\frac{l_i}2}.$$ Thus $r_i$ converges to $\infty$ uniformly.
\end{proof}

\begin{lemma}\label{d} The following holds for some positive
finite numbers $f_1,f_2,f_3,f_4,f_5$:

(1) if $(l_1,l_2,l_3)$ converges to $(\infty, f_1, f_2)$, then
$(\theta_1,\theta_2,\theta_3)$ converges to $(\infty, f_3, f_4)$;

(2) if $(l_1,l_2,l_3)$ converges to $(\infty,\infty, f_5)$, then
$\theta_3$ converges to 0;

(3) if $(l_1,l_2,l_3)$ converges to $(\infty,\infty,\infty)$, then
we can chose a subsequence of $(l_1,l_2,l_3)$ such that at least
two of $\theta_1,\theta_2$ and $\theta_3$ converge to 0.
\end{lemma}

\begin{proof}
(1) By the cosine law we have  $$\cosh\theta_1 =\frac{ \cosh l_1 +
\cosh l_2\cosh l_3}{\sinh l_2 \sinh l_3},$$ if $\lim(l_1,l_2,l_3)=
(\infty, f_1, f_2)$, we have $\lim\cosh\theta_1 = \infty,$ or
$\lim\theta_1 = \infty.$ And since $\lim\frac{\cosh l_1}{\sinh
l_1}=1,$
$$\lim\cosh \theta_2 =\lim\frac{ \cosh l_2 +
\cosh l_1\cosh l_3}{\sinh l_1 \sinh l_3}=\frac{\cosh f_2}{\sinh
f_2}>1.$$ Thus $\lim\theta_2$ is a positive finite number. By
symmetry $\lim\theta_3$ is a positive finite number.

(2) If $\lim(l_1,l_2,l_3)=(\infty,\infty, f_5)$, we have
$$\lim\cosh\theta_3 =\lim\frac{ \cosh l_3 +
\cosh l_1\cosh l_2}{\sinh l_1 \sinh l_2}=\lim\frac{\cosh
l_3}{\sinh l_1 \sinh l_2}+1=1.$$ Thus $\lim\theta_3=0$.

(3) If $\lim(l_1,l_2,l_3)=(\infty,\infty,\infty)$, we have
$$\lim\cosh\theta_i=\lim\frac{ \cosh l_i + \cosh l_j\cosh l_k}{\sinh l_j
\sinh l_k}=\lim\frac{\cosh l_i}{\sinh l_j \sinh l_k}+1$$
$$=\lim\frac{2e^{l_i}}{e^{l_j+l_k}}+1=\lim2e^{l_i-l_j-l_k}+1.$$
Since $\lim e^{l_i-l_j-l_k}e^{l_j-l_i-l_k}=\lim e^{-2l_k}=0,$ by
taking subsequence of $(l_1,l_2,l_3)$, we may assume $\lim
e^{l_i-l_j-l_k}$ and $\lim e^{l_j-l_i-l_k} $ exist. Then one of
$\lim e^{l_i-l_j-l_k}$ and $\lim e^{l_j-l_i-l_k} $ is 0. Hence at
least two of $\lim\theta_1, \lim\theta_2$ and $\lim\theta_3$ are 0.
\end{proof}

\section{Proof of Theorem \ref{negative}}

\begin{lemma}\label{sum} If $a > 0$, then for any real number x,
we have $$\int_0^{a+x} \cosh^{\lambda}(t) dt + \int_0^{a-x}
\cosh^{\lambda}(t)dt > 0.$$
\end{lemma}

\begin{proof} Let $f(a)$ be the function of the left hand side of the inequality. We see
$f'(a)=\cosh^{\lambda}(a+x)+\cosh^{\lambda}(a-x)>0.$ And $f(0)=0.$
Hence $f(a)>0$ for $a>0.$
\end{proof}

\begin{proof}[Proof of Theorem \ref{negative}] We denote the polytope defined by the
inequalities in condition 1, 2, 3 by $P_{\lambda}.$ First we claim
$\Psi_{\lambda}(T(S)) \subset P_{\lambda}$. Indeed, fix a
hyperbolic metric $l \in T(S)$. For any edge $e,$ let $r,r'$ be
the r-coordinates of A-arcs facing $e$, then
$$\psi_\lambda(e)=\int_0^{r} \cosh^{\lambda}(t) dt + \int_0^{r'}
\cosh^{\lambda}(t)dt < 2 \int_{0}^{\infty} \cosh^{\lambda}(t)
dt.$$ Thus the condition 1 holds.

Given an edge path $(H_0, e_1, H_1,..., e_n, H_n)$, for
$i=1,...,n-1,$ let $a_i$ be the length of the A-arc in $H_i$
adjacent to $e_i$ and $e_{i+1}$. Denote the lengths of A-arcs in
$H_i$ facing $e_i$ and $e_{i+1}$ by $b_i$ and $c_i$ respectively
as labelled in Figure \ref{path} (a).

\begin{figure}[htbp]
\begin{center}
\includegraphics[scale=.55]{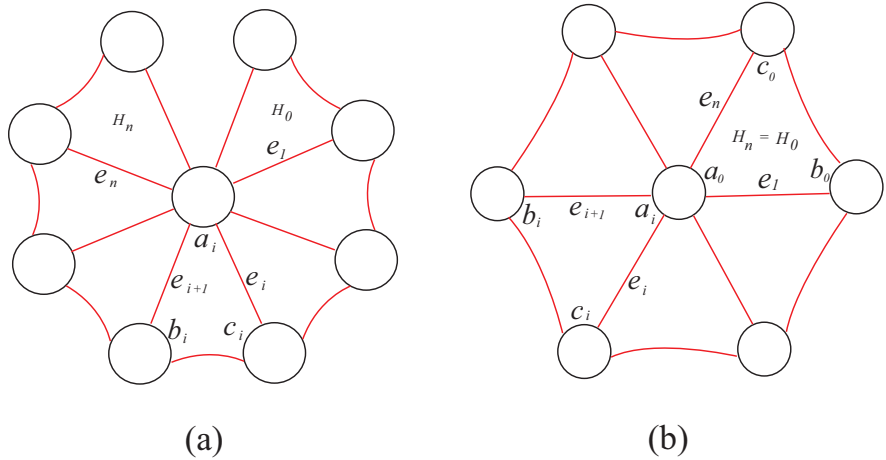}
\end{center}
\caption{\label{path}(a) An example of an edge path with lengths of A-arcs labeled. (b) An example
of an edge cycle with lengths of A-arcs labeled.}
\end{figure}

Then by definition
$$\psi_\lambda(e_1)=\int_0^{r} \cosh^{\lambda}(t) dt +
\int_0^{\frac{a_1+c_1 -b_1}{2}} \cosh^{\lambda}(t)dt,$$ where $r$
is the r-coordinate of the A-arc in $H_0$ facing $e_1$. For
$i=2,...,n-1,$
$$\psi_\lambda(e_i)=\int_0^{\frac{a_{i-1}+b_{i-1}-c_{i-1}}{2}} \cosh^{\lambda}(t) dt +
\int_0^{\frac{a_i+c_i-b_i}{2}} \cosh^{\lambda}(t)dt.$$ And
$$\psi_\lambda(e_n)=\int_0^{\frac{a_{n-1}+b_{n-1}-c_{n-1}}{2}} \cosh^{\lambda}(t) dt +
\int_0^{r'} \cosh^{\lambda}(t)dt$$ where $r'$ is the r-coordinate
of the A-arc in $H_n$ facing $e_n.$

Hence by Lemma \ref{sum},
\begin{eqnarray*}
&&\sum_{i=1}^n \psi_\lambda(e_i)= 
\sum_{i=1}^{n-1}(\int_0^{\frac{a_i+c_i-b_i}{2}}
\cosh^{\lambda}(t)dt 
+ \int_0^{\frac{a_i+b_i-c_i}{2}}
\cosh^{\lambda}(t)dt) \\ &+&\int_0^{r} \cosh^{\lambda}(t) dt
+\int_0^{r'} \cosh^{\lambda}(t)dt \\
&>&\int_0^{r}\cosh^{\lambda}(t) dt + \int_0^{r'}
\cosh^{\lambda}(t)dt\\
&>& 2 \int_{0}^{-\infty} \cosh^{\lambda}(t) dt =
-2\int_{0}^{\infty} \cosh^{\lambda}(t) dt.
\end{eqnarray*}
Thus the condition 2 holds.

Furthermore, if $(H_0, e_1, H_1,..., e_n, H_n=H_0)$ is an edge
cycle, $H_0$ contains both $e_1$ and $e_n$. Let $a_0$ be the
length of A-arc in $H_0$ adjacent to $e_1$ and $e_n$, $b_0,c_0$ be
the lengths of A-arcs facing $e_n$ and $e_0$ respectively as
labelled in Figure \ref{path} (b). Thus the r-coordinates are
$r=\frac{a_0+b_0-c_0}2$ and $r'=\frac{a_0+c_0-b_0}2.$ Hence
$$\sum_{i=1}^n \psi_\lambda(e_i)>\int_0^{r}\cosh^{\lambda}(t) dt +
\int_0^{r'} \cosh^{\lambda}(t)dt>0$$ by Lemma \ref{sum}. Thus the
condition 3 holds.

Now by Theorem \ref{embedding}, $\Psi_{\lambda}:T(S) \to
P_{\lambda}$ is an embedding. Therefore $\Psi_{\lambda}(T(S))$ is
open in $P_{\lambda}.$ We only need to show it is also closed in
$P_{\lambda}.$ This will finish the proof since $P_{\lambda}$ is
connected.

Take a sequence $l^{(m)} \in T(S)$ so that $\lim_{m \to
\infty}\Psi_{\lambda}(l^{(m)})=z\in P_{\lambda}.$ By taking
subsequence, we may assume that $\lim_{m \to \infty} l^{(m)}\in
[0, \infty]^E$ exists and the length of each A-arc converges into
$[0, \infty]$. We only need to show that $\lim_{m \to \infty}
l^{(m)} \in (0, \infty)^E=T(S).$ This will finish the proof since
$z=\Psi_{\lambda}(\lim_{m \to \infty} l^{(m)}).$

Suppose otherwise that there is an edge $e \in E$ so that $
\lim_{m \to \infty} l^{(m)}(e) \in \{0, \infty\}$. We will discuss
two cases.

Case 1, $\lim_{m \to \infty} l^{(m)}(e) =0$ for some $e \in E$.
Let $H,H'$ be the hexagons sharing $e$ and $r^{(m)},r'^{(m)}$ be
the r-coordinates of the A-arcs in $H,H'$ facing $e$. Then by
Lemma \ref{r}, $\lim_{ m \to \infty} r^{(m)}\rightarrow \infty,
\lim_{ m \to \infty} r'^{(m)}\rightarrow \infty.$ Then $$z(e)=
\lim_{ m \to \infty} (\int_{0}^{r^{(m)}} \cosh^{\lambda}(t)
dt+\int_{0}^{r'^{(m)}} \cosh^{\lambda}(t) dt ) = 2
\int_{0}^{\infty} \cosh^{\lambda}(t) dt.$$ This is impossible
since $z\in P_{\lambda}$ must satisfy the condition 1.

Due to case 1, we can assume $\lim_{m \to \infty}
l^{(m)}\in(0,\infty]^E.$

Case 2, $\lim_{m \to \infty} l^{(m)}(e) =\infty$ for some $e \in
E$. Define the subset $E_{\infty}=\{e\in E | \lim_{m \to \infty}
l^{(m)}(e)=\infty\}.$ We construct a graph $G$ as follows. A
vertex of $G$ is a hexagon with at least one edge in $E_{\infty}$.
There is a \it dual-edge \rm in $G$ joining two vertexes if and
only if the two hexagons corresponding to the vertexes share an
edge in $E_{\infty}$. The degree of a vertex of the graph $G$
can only be 1, 2 or 3. Actually a vertex of degree 1, 2 or 3 is
corresponding to the hexagon of type (1),(2) or (3) in Lemma
\ref{d} respectively.

We smooth the graph $G$ at vertexes as follows. At a vertex of
degree 1, we replace the small neighborhood of
the vertex in $G$ by a short smooth curve tangent to the unique
dual-edge incident to the vertex as in Figure \ref{modify} (a). At a
vertex $v$ of degree of 2 or 3, every two dual-edges
$\overline{e}_1,\overline{e}_2$ incident to $v$ correspond to two
edges $e_1,e_2$ in a hexagon. If the length of the A-arc adjacent
to $e_1,e_2$ converges to 0, we replace the small neighborhood of
the vertex $v$ in $G$ by a short smooth curve tangent to
$\overline{e}_1,\overline{e}_2$. According to Lemma \ref{d}, every
vertex of degree 2 can be smoothed as in Figure \ref{modify} (b)
and there are two cases for a vertex of degree 3 according to the
lengths of 2 or 3 A-arcs converge to 0 as in Figure \ref{modify}
(c).

\begin{figure}[htbp]
\begin{center}
\includegraphics[scale=.6]{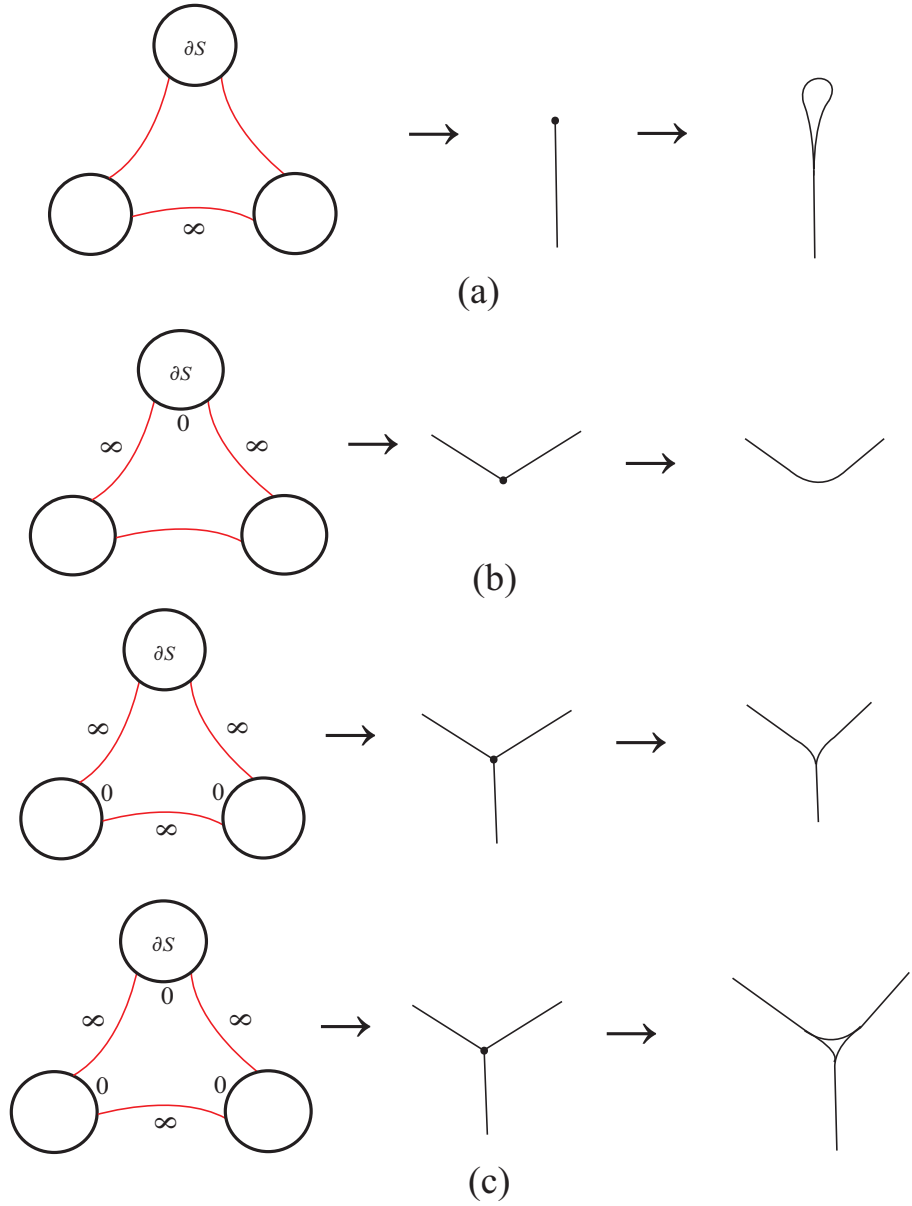}
\end{center}
\caption{\label{modify}Smooth graph G at a vertex of degree 1 (a), degree 2 (b), degree 3 (c).}
\end{figure}

We denote by $G'$ the graph smoothed at vertexes and the
dual-edges of $G'$ are the dual-edges of $G$. We claim that there
exists a smooth closed curve in $G'$ such that every dual-edge repeats at
most twice in the closed curve. In fact, we give every dual-edge of $G'$
an arbitrary orientation. Pick up any smooth closed curve in $G'$ which may
contains infinite number of dual-edges. If there exists an
dual-edge $\overline{e}$ repeats with the same orientation in the
closed curve, there is another smooth closed curve starting and ending at $\overline{e}$.
By this procedure we can reduce the number of dual-edges of a
closed curve. At last we obtain a smooth closed curve in $G'$ such that every dual-edge
repeats at most twice.

This smooth closed curve in $G'$ corresponds a fundamental edge path or
fundamental edge cycle in the ideal triangulation. First assume it
is a fundamental edge path $(H_0, e_1, H_1,..., e_n, H_n).$ Since
the degree of the vertex corresponding to $H_0$ (or $H_n$) is 1,
the lengths of other two edges other than $e_1$ (or $e_n$) converge
to positive finite numbers in the sequence of metric $l^{(m)}$. By
Lemma \ref{d}(1) the r-coordinate of the A-arc in $H_0$ (or $H_n$)
facing $e_1$ (or $e_n$) converges to $-\infty.$ By the
construction of the edge path, the length of A-arc adjacent to
$e_i$ and $e_{i+1}$ converges to 0 for $i=1,...,n-1.$ And we
denote $b_i,c_i$ the limit of lengths of A-arcs in $H_i$ facing
$e_i$, $e_{i+1}$ respectively, see Figure \ref{path}(a).

Hence $$z(e_1)=\int_0^{-\infty} \cosh^{\lambda}(t) dt +
\int_0^{\frac{c_1 -b_1}{2}} \cosh^{\lambda}(t)dt.$$ For
$i=2,...,n-1,$
$$z(e_i)=\int_0^{\frac{b_{i-1}-c_{i-1}}{2}} \cosh^{\lambda}(t) dt +
\int_0^{\frac{c_i-b_i}{2}} \cosh^{\lambda}(t)dt.$$ And
$$z(e_n)=\int_0^{\frac{b_{n-1}-c_{n-1}}{2}} \cosh^{\lambda}(t) dt +
\int_0^{-\infty} \cosh^{\lambda}(t)dt.$$

Hence $$\sum_{i=1}^n z(e_i)=2\int_0^{-\infty} \cosh^{\lambda}(t)
dt +\sum_{i=1}^{n-1}(\int_0^{\frac{c_i-b_i}{2}}
\cosh^{\lambda}(t)dt + \int_0^{\frac{b_i-c_i}{2}}
\cosh^{\lambda}(t)dt) $$
$$= -2 \int_{0}^{\infty}\cosh^{\lambda}(t) dt.$$ This is impossible since $z\in P_{\lambda}$ must
satisfy the condition 2.

If the smooth closed curve in $G'$ corresponds to a fundamental edge cycle $(H_0,
e_1, H_1,..., e_n, H_n=H_0)$, the length of A-arc in $H_0$
adjacent to $e_1$ and $e_n$ is 0. Denote $b_0,c_0$ the lengths of
A-arcs facing $e_n$ and $e_0$, see Figure \ref{path}(b). Thus
$$z(e_1)=\int_0^{\frac{b_0-c_0}{2}} \cosh^{\lambda}(t) dt +
\int_0^{\frac{c_1 -b_1}{2}} \cosh^{\lambda}(t)dt,$$ 
$$z(e_n)=\int_0^{\frac{b_{n-1}-c_{n-1}}{2}} \cosh^{\lambda}(t) dt +
\int_0^{\frac{c_0 -b_0}{2} } \cosh^{\lambda}(t)dt.$$ And as in the case of 
fundamental edge path, for $i=2,...,n-1,$
$$z(e_i)=\int_0^{\frac{b_{i-1}-c_{i-1}}{2}} \cosh^{\lambda}(t) dt +
\int_0^{\frac{c_i-b_i}{2}} \cosh^{\lambda}(t)dt.$$ Hence
$\sum_{i=1}^n z(e_i)=0$. This is impossible since $z\in
P_{\lambda}$ must satisfy the condition 3.

We finish the proof of $\Psi_{\lambda}(T(S))=P_{\lambda}$. 
Since there are only finite many fundament edge paths or fundament
edge cycles in an ideal triangulation, $\Psi_{\lambda}(T(S))$ is
defined by finite many inequalities in condition 1, 2, 3. Thus it
is a open convex polytope.

The statement $\Psi_{\lambda_1}(T(S))\subset
\Psi_{\lambda_2}(T(S))\subset\Psi_0(T(S))=\Psi_{\lambda_3}(T(S))$
for $\lambda_1<\lambda_2<0<\lambda_3$ is obvious since the
function $\int_{0}^{\infty} \cosh^{\lambda}(t) dt$ is increasing
in $\lambda$ and it is $\infty$ when $\lambda\geq 0.$

Since $0<\cosh^{\lambda}(t)<\cosh^{-1}(t)$ for $\lambda<-1$ and
$\int_{0}^{\infty} \cosh^{-1}(t) dt<\infty,$ by Lebesgue's
dominated convergence theorem, we have
$$\lim_{\lambda\to-\infty}\int_{0}^{\infty} \cosh^{\lambda}(t) dt=
\int_{0}^{\infty}\lim_{\lambda\to-\infty} \cosh^{\lambda}(t)
dt=0.$$ Thus the intersection
$\cap_{\lambda=0}^{-\infty}\Psi_{\lambda}(T(S))$ is the set of points $z\in \mathbf{R}^E$
satisfying $z(e) < 0$ for each
edge $e$ and $\sum_{i=1}^n z(e_i) >0$ for each fundamental edge path $(H_0, e_1, H_1,...,$ $e_n,
H_n).$ It is an empty set.
\end{proof}

\section*{Acknowledgement} The author would like to thank his advisor, Feng
Luo, for suggesting this problem and helpful discussions. 
This work is partially supported by NSF Grant \#0604352.

\end{document}